\documentclass[thmsb,11pt,final,notitlepage,oneside,letterpaper]{article}
\usepackage{amssymb}

\usepackage{sw20lart}

\input{tcilatex}
\begin{document}

\begin{center}
{\Large \textbf{Minimum aberration} }$2^{k-p}${\Large \ \textbf{designs of
resolution }III} \vspace{0.5cm}

\begin{tabular}{cc}
Jes\'{u}s Juan & J. Gabriel Palomo \\ 
Lab. de Estad\'{i}stica (E.T.S.I. Industriales) & Dep. Matem\'{a}tica
Aplicada A.T. \\ 
Univ.\ Polit\'{e}cnica de Madrid, Spain & Univ.\ Polit\'{e}cnica de Madrid,
Spain \\ 
E-mail: \texttt{jjuan@etsii.upm.es} & E-mail: \texttt{palomo.san@euatm.upm.es%
}
\end{tabular}
\bigskip

ABSTRACT \bigskip

\begin{minipage}{12cm}
\small
In this article we prove several important properties of $2^{k-p}$  minimum aberration (MA) 
designs with $k \ge  n/2$,  where $n=2^{k-p}$ is 
the number of runs. We develop a simple method to build MA designs of resolution III. 
Furthermore, we  introduce a 
simple relationship, based on product of polynomials, for computing 
their word-length patterns. 
 \smallskip


\end{minipage}
\end{center}

\section{Introduction}

$2^{k-p}$ designs are commonly used in experimentation, where $k$ is the
number of factors and $n=2^{k-p}$ is the number of runs or observations.
These designs are economical, statistically efficient and simple to analyze
(see Box and Hunter, 1961).

In many applications, especially in industrial experiments, it is often
necessary to determine which factors among a large number of candidates
could affect a particular response. For a given number of factors $k$, a
usual procedure is to select a first fractional $2^{k-p}$ design to
accommodate the factors with the smallest number of runs. Confusion of
effects may lead to more than one plausible explanation of data. The initial
design should be chosen in such a way that simplifies as much as possible
the analysis stage and reduces the size of subsequent experimentation; see
Box, Hunter and Hunter (1978) or Tiao (2000) for examples. When the
experimenter has little knowledge about the relative sizes of the factorial
effects, the minimum aberration (MA)\ criterion selects designs with good
overall properties. A detailed discussion on the MA criterion can be found
in Fries and Hunter (1980) and Chen, Sun and Wu (1993). Since Fries and
Hunter (1980) introduced this criterion, many articles have been devoted to
find MA designs and to study their characterization and structure, for
example, Franklin (1984), Chen and Wu (1991), Chen (1992) and Tang and Wu
(1996).

In addition to their direct practical application, minimum aberration
designs are also useful in the construction of more complex designs. For
example, Ankenman (1999) includes four-level factors in a design, and
Bingham and Sitter (1999) and Huang, Chen and Voelkel (1998) use them to
obtain minimum aberration two-level split-plot designs.

MA designs can be obtained with reduced computational effort, particularly
when the number of runs is low (16, 32 or 64). In this paper we propose an
original procedure to obtain MA designs of resolution III in a simple way.
This result could be useful from a practical point of view when the number
of runs is high. That is not, however, our main interest. Following the work
of Tang and Wu (1996), we wish to understand the structure of MA designs,
establishing some of their characteristic properties, and to place them
within the global structure of designs with a fixed number of runs. We will
see graphically that this structure presents interesting properties.

In this paper we shall only consider screening designs. Thus, we will
consider $2^{k-p}$ designs which allow us to study $k$ factors with a
minimum number of runs $n=2^{k-p},$ i.e., $n/2\leq k\leq n-1.$ There are two
reasons for this choice: first, they are preferred when the number of
factors to be analyzed, $k$, is large; this is precisely the case when the
search for a good design is more complex. Second, using a symmetry property
(Tang and Wu , 1996) the properties of designs with $k<n/2$ can be deduced
from those with $k\geq n/2$. As a consequence, a complete characterization
of designs having $k\geq n/2$ provides a characterization for all designs
with $n$ runs. Moreover, all designs with $k\geq n/2$, except for one, are
of resolution $III$ (for $k=n/2),$ and their comparison would require
precise criteria to be able to discriminate among them.

Finally, and in addition to providing a characterization of MA designs for $%
k\geq n/2,$ we also provide a way to replace the tedious computation of the 
\emph{word-length pattern} vector by the computation of an elementary
product of polynomials.

This paper is organized in four sections. Section 2 introduces the problem
and the notation. In Section 3 the main result for this paper is presented,
namely the necessary condition for any $2^{k-p}$ design with $k\geq n/2$ to
be of minimum aberration. Sections 4 and 5 are concerned with the
computation of the \emph{word-length pattern}, and from the results in these
sections we complete the conditions which characterize MA designs.

\section{Definitions and motivation}

A complete two-level factorial design with $m$ factors $%
B_{1},B_{2},...,B_{m} $ contains the $n=2^{m}$ possible combinations of
factor levels. Usually, - and + signs are employed to indicate the two
levels of each factor and the design is presented in a \emph{design matrix},
which has $m$ columns, one per factor, and $2^{m}$ rows corresponding to the
factor treatments. The $m$ columns are independent in that none of them can
be obtained as a product of the remaining ones. To simplify the
presentation, we shall call these columns or factors the \emph{basic} ones.
By multiplying these columns in all possible ways we can generate up to $%
s=2^{m}-m-1$ additional ones, which will be identified by the combination
(juxtaposition) of the basic factors which generate them $\pi
_{i}(B_{1},B_{2},...,B_{m})$ or more concisely as $\pi _{i}$. The individual
elements of column $\pi _{i}$ are obtained by multiplying the corresponding
elements of the columns specified in $\pi _{i}(B_{1},B_{2},...,B_{m}).$ The
collection of the $2^{m}-1$ columns 
\[
H_{m}=\{B_{1},B_{2},...,B_{m},\pi _{1},\pi _{2},...,\pi _{s}\} 
\]
is the key set through which any $2^{k-p}$ designs with $2^{m}$ runs can be
defined. In coding theory the set $H_{m}$ is known as a Hamming code. For
example, if we use $A,B,C,D$ for the four initial columns of the complete $%
2^{4}$ design, we obtain 
\[
H_{4}=\{A,B,C,D,AB,AC,AD,BC,BD,CD,ABC,ABD,ACD,BCD,ABCD\}. 
\]

Multiplying any column of $H_{m}$ by itself we obtain the $I$ column, where
all signs are +. The set $\left\{ I\right\} \cup H_{m}$ with the defined
operation has a structure of an abelian group, where $I$ is the null
element. The result of the product of sign columns can also be derived by
looking at the products of the combinations of letters which represent them.
For example, the product of columns $ABC$ and $ABD$ produces the column $CD,$
since $AA=I$, $BB=I.$

\begin{definition}
\label{Def_1} A subset $d$ of $k$ columns of $H_{m}$ represents a $2^{k-p}$
design.
\end{definition}

For example, the following subset of $9$ columns of $H_{4}$ 
\begin{equation}
d_{1}=\left\{ A,B,C,D,ABC,ABD,ACD,BCD,ABCD\right\}  \label{dis1}
\end{equation}
is a $2^{9-5}$ design$.$

The advantage of this definition derives from its generality; it enables us
to define any design with $n=2^{m}$ runs, whether it is a complete factorial
design or a fractional design $2^{k-p}$ ($k>m)$. Moreover, through this
definition we can make use the tools of set theory. We will then say that a
given design contains another, is the union of two others, and so forth.

Many subsets formed by different elements of $H_{m}$ provide equivalent
designs which for our work will be considered the same. Two designs are
equivalent if it is possible to define an isomorphism between them (see for
example, Tang and Wu, 1996). Thus, Pu (1989) proved that there are only five
different $2^{9-5}$ designs, so that any subset of $H_{4}$ with nine
elements is equivalent or isomorphic to one of those designs. If we consider
all designs, with any number of factors, the number of different designs
with $n=16$ runs is reduced to $45$ ($46$ if we include $\varnothing $). In
figure 1 we present a graph showing the relationship between them. Each
point or node in the graph is a design. The starting point is the empty set $%
\varnothing $ ($k=0$); for $k=1$ all designs are equivalent and are
represented by a single node, and the same occurs for $k=2.$ For $k=3$ there
are two different designs, $\{A,B,C\}$ and $\{A,B,AB\};$ for $k=4$ there are
three, and so forth up to $k=15$ which is the $H_{4}$ design. The lines
connect pairs of designs which only differ in one element, which in some
cases (solid lines) is indicated on the arc. The design with $k+1$ columns
is obtained by adding the column which appears on the arc that connects it
to its predecessor. The designs for $k>8$ can be obtained by symmetry,
starting from $H_{4}$ and eliminating the columns specified in the symmetric
arcs.

One of the most important properties of this graph is its symmetry, due to
the fact that if two designs $d_{1}$ and $d_{2}$ are equivalent, so are
their complementaries $\bar{d}_{1}$ and $\bar{d}_{2}$. The notation $\bar{d}$
will be used to represent the complementary of a given design $d$, where $%
\bar{d}=H_{m}\backslash d$. Thus, within the subsets of $H_{4}$, there are
four different designs for $k=5$, there exist four others for $k=10$, which
correspond to the complementaries of the first ones. In figure 1 the
complementary of a design with $k$ columns is the symmetric with $15-k$
columns. The relationships between complementary designs are analyzed by
Tang and Wu (1996). In accordance with this, the properties of one half of
these designs ($k<n/2)$ can be deduced from the other half $k\geq n/2$. In
this paper we shall concentrate in the analysis of MA designs for $k\geq n/2$%

The goodness of a design is a consequence of the relations of dependence
between its columns. Let us consider the design $d\subset H_{m},$ with
columns $d=\{\xi_{1},\xi _{2},...,\xi _{k}\}.$ A word of length $i$ consists
of $i$ elements $\xi _{j_{1}},\xi _{j_{2}},...,\xi _{j_{i}}$ from $d$ such
that 
\begin{equation}
\xi _{j_{1}}\xi _{j_{2}}\cdots \xi _{j_{i}}=I  \label{word}
\end{equation}
where $I$ denotes the column having all $+$'s. The set of all distinct words
formed by products involving elements of $d$ gives the \emph{defining
relation} of the design. The words in the defining relation correspond to
all the interactions of the $k$ factors that are confounded with the mean
(represented by column $I$). The vector 
\[
W(d)=[a_{1}(d),a_{2}(d),...,a_{k}(d)] 
\]
is called the \emph{word-length pattern} (WLP) of the design $d$, where $%
a_{i}(d)$ is the number of words of length $i$. We shall use in the future $%
a_{i}$ instead of $a_{i}(d)$ whenever this does not lead to confusion.

The resolution of a design is the length of its shortest word. Maximum
resolution is the usual criterion to select a design. However, for given $k$
and $p$, and specially when these values are large, there may exist many
different maximum resolution $2^{k-p}$ designs. For instance, the five
different $2^{9-5}$ designs are of resolution $III$. Additionally, for $%
k\geq n/2$ all $H_{m}$ designs but one (the exception corresponds to $k=n/2$
and has resolution \emph{$IV$}) are of resolution $III$.

Fries and Hunter (1980) introduced the notion of \emph{aberration} as a more
powerful method to compare fractional factorial designs. For two $2^{k-p}$
designs $d_{1}$ and $d_{2},$ suppose $r$ is the smallest value such that $%
a_{r}(d_{1})\neq a_{r}(d_{2})$. We say that $d_{1}$ has smaller aberration
than $d_{2}$ if $a_{r}(d_{1})<a_{r}(d_{2}).$ The design $d$ is of minimum
aberration (MA) if there is no other design with the same number of factors
and runs having less aberration than $d.$

This previous analysis is of practical interest for designs with $k>m,$
which are called fractions. The usual way of defining a fraction is the
following: first $m$ initial factors $b=\{B_{1},B_{2},...,B_{m}\}$ are taken
to form the $2^{m}$ design, and then $p$ additional factors $%
B_{m+1},B_{m+2},...,B_{m+p}$ are assigned to a subset $\left\{ \xi
_{1}^{\prime },\xi _{2}^{\prime },...,\xi _{p}^{\prime }\right\} ,$ selected
from $\bar{b}=H_{m}\backslash b.$ The selection is denoted by 
\begin{equation}
\left. 
\begin{array}{c}
B_{m+1}=\xi _{1}^{\prime }(B_{1},B_{2},...,B_{m}) \\ 
B_{m+2}=\xi _{2}^{\prime }(B_{1},B_{2},...,B_{m}) \\ 
\vdots \\ 
B_{m+p}=\xi _{p}^{\prime }(B_{1},B_{2},...,B_{m})
\end{array}
\right\}  \label{canonica}
\end{equation}
The design $d$ is thus 
\begin{equation}
d=\left\{ B_{1},B_{2},...,B_{m},\xi _{1}^{\prime },\xi _{2}^{\prime
},...,\xi _{p}^{\prime }\right\}  \label{d}
\end{equation}
This definition requires the inclusion of the $m$ basic columns. Any set of $%
H_{m}$ containing $m$ independent columns provides designs equivalent to
those defined from the above criterion. If we assume that this method has
been applied to define the $2^{k-p}$ design, the fraction is uniquely
determined for $\xi _{1}^{\prime },\xi _{2}^{\prime },...,\xi _{p}^{\prime }$%
. Multiplying each relation $B_{m+j}=\xi _{j}^{\prime }$ by $B_{m+j}$ we
obtain $I=\xi _{i}^{\prime }B_{m+j}$. The $p$ terms $\xi _{i}^{\prime
}B_{m+j}$ are called generators of the fraction and represent effects
confounded with the mean. These relationships are usually rewritten as 
\begin{equation}
I=\xi _{1}^{\prime }B_{m+1}=\xi _{2}^{\prime }B_{m+2}=\cdots =\xi
_{p}^{\prime }B_{m+p}.  \label{generadores}
\end{equation}

The product of two generators of (\ref{generadores}) provides another effect
that is also confounded with the mean. The set of distinct words formed by
all possible products involving the $p$ generators gives the \emph{defining
relation} of the fraction. It is usual to refer to letters instead of
factors or columns. A letter is any of the labels $B_{i}$ used to denote a
factor. In this setting, as before, the product of two identical letters is
the identity. The \emph{word length} is simply the number of letters of a
word.

For instance, if we call $E,F,G,H$ and $J$ the five additional factors, the
generators of the design (\ref{dis1}) are 
\begin{equation}
I=ABC\underline{E}=ABD\underline{F}=ACD\underline{G}=BCD\underline{H}=ABCD%
\underline{J}.  \label{ej1}
\end{equation}
(The non basic factors have been underlined). By taking all possible
products in which the five generators of (\ref{ej1}) are present, we obtain
the 32 words (31 if we exclude $I$) which form the defining relation of the
fraction. From the lengths of the words we obtain the following WLP, 
\begin{equation}
(0,0,4,14,8,0,4,1,0).  \label{wlpej1}
\end{equation}
The resolution of this design is $III$, which is the maximum for $2^{9-5}$
designs.

For the cases to be considered, $k\geq n/2$, any subset has always $m$
independent columns and always leads to designs which can be defined by
means of the scheme described in (\ref{canonica}). In the following sections
we shall use any of the two forms of identifying a design.

The WLP is an essential instrument to evaluate a $2^{k-p}$ design and its
computation requires $p$ independent generators. In section 5 we shall
encounter designs which are defined for subsets of $H_{m}$ that do not
include the basic columns such as 
\[
d_{2}=\{A,BC,BD,CD,ABC,ABD,ACD,BCD,ABCD\}. 
\]
A simple method to obtain the generators consists in transforming $d_{2}$ by
means of an isomorphism into another set with structure (\ref{d}). For
instance, the isomorphism $\mathcal{T}$ defined by $\mathcal{T(}A)=A,%
\mathcal{T(}B)=BCD,\mathcal{T(}C)=ACD$ and $\mathcal{T(}D)=ABD$ transforms $%
d_{2}$ into the equivalent design 
\[
\{A,B,C,D,AB,AC,AD,BC,ABC\}. 
\]
Using $E,F,G$ and $J$ as labels for the non basic factors, the generators
are 
\[
I=AB\underline{E}=AC\underline{F}=AD\underline{G}=BC\underline{H}=ABC%
\underline{J}. 
\]
The WLP of $d_{2}$ is $(0,0,8,10,4,4,4,1,0)$. Since $a_{3}(d_{1})=4$ and $%
a_{3}(d_{2})=8,$ $d_{1}$ is preferred from the aberration point of view.

Tang and Wu (1996) propose a new approach to characterize MA designs in
terms of their complementary designs. Let $H_{m}=d\cup \bar{d}$ $,$ they
argue that when the elements in $\bar{d}$ are more ``dependent'', those in $%
d $ should be less ``dependent'' and thus may have less aberration. They
developed a general theory to support this intuition. Here we give a
different version of Tang and Wu's intuition based in the concept of \emph{%
rank} which we define as follows.

\begin{definition}
The rank of a set $f$ of columns of $H_{m}$ is the number of independent
columns it contains.
\end{definition}

The rank of $H_{m}$ is $m$ and hence the rank of any subset $f\subset H_{m}$
is never larger than $m.$ Moreover, the maximum number of elements of a
subset $f$ with rank $v$ is $2^{v}-1.$ Therefore, if the number of elements
of a set is $h$, its minimum rank is 
\[
v_{h}=\lceil \log _{2}(h+1)\rceil , 
\]
where $\lceil (a)\rceil $ is the smallest integer larger or equal to $a.$

The maximum resolution (and also the MA) fractions always have maximum rank
equal to $m$. In the next section we shall also see that the complementary
design of a MA design has minimum rank.

Let us go back to figure 1 to identify the MA designs. For each value of $k$
in the figure, we have identified with a black node the MA design. The
complementary (symme\-tric) of a MA design is the worse design from the
aberration point of view. All designs with $k\geq 8$ are of resolution $III$
with the exception of a design for $k=8$ having resolution IV (the circled
one). This design is called \emph{saturated design with resolution IV} and
plays an essential role when building minimum aberration designs. We can
observe that all minimum aberration designs for $k>8$ start from it. This
property and other related ones are generalized in the following sections.

\section{Screening designs $(k\geq n/2)$ of minimum aberration}

Let $d$ be the $2^{k-p}$ design with $n/2\leq k\leq n-1,$ defined by the
generators $I=\xi _{1}^{\prime }B_{m+1}=\xi _{2}^{\prime }B_{m+2}=\cdots
=\xi _{p}^{\prime }B_{m+p}.$ Our goal is to choose $\xi _{1}^{\prime },\xi
_{2}^{\prime },...,\xi _{p}^{\prime }$ in such a way that the design $d$ is
of minimum aberration.

The solution for the case $k=n/2$ is well known. Given the basic columns $%
B_{1},B_{2},...,B_{m},$ it consists in assigning the $(n/2-m)$ non-basic
factors to the columns obtained as the product of an odd number of basic
ones. There is only one possible design of resolution IV and hence it is of
minimum aberration.

Let $\varphi _{1},\varphi _{2},...,\varphi _{q}$ with $q=2^{m-1}-m,$ be all
the odd combinations of the basic columns $B_{1},B_{2},...,B_{m}.$ The
generators of the $2_{IV}^{l-q}$ design with $l=2^{m-1}$ are 
\begin{equation}
I=\varphi _{1}B_{m+1}=\varphi _{2}B_{m+2}=\cdots =\varphi _{q}B_{m+q}
\label{saturado}
\end{equation}
From now on we will denote this design as $O_{m}.$

The extreme case is obtained when $k=n-1,$ and corresponds to the maximum
number of factors which can be analyzed with a regular fraction of $n$ runs.
This design is called a saturated design of resolution $III$ (or just a
saturated design) and will be denoted as $H_{m},$ since it contains all its
columns.

When the number of factors to be studied is $n/2<k<n-1,$ the maximum
resolution of the design is $III$. The resolution criterion does not
discriminate among designs with this number of factors and it is necessary
to use the aberration criterion to differentiate between them. The following
theorem gives a condition for a design with a number of factors $k>n/2$ to
be of MA.

\begin{theorem}
Let $d\subset H_{m}$ be a $2^{k-p}$ design and $\bar{d}$ its complementary, $%
H_{m}=d\cup \bar{d}$. A necessary condition for $d$ to be of minimum
aberration is that the set $\bar{d}$ has minimum \emph{rank.}
\end{theorem}

The proof, based on a combinatorial argument, is given in appendix 1.

It should be noted that theorem 1 does not limit the number $h$ of columns
of $\bar{d}$. When $h>n/2$, the rank of the set $\bar{d}$ of columns to be
eliminated is always $m$ and the theorem is trivial. The result is important
for $h<n/2$, which corresponds to designs of resolution $III$. We wish to
emphasize the interest of this result when looking for MA fractions. For
example, the $2^{12-8}$ design with generators 
\begin{equation}
I=ABC\underline{E}=ABD\underline{F}=ACD\underline{G}=BCD\underline{H}=AD%
\underline{J}=BD\underline{K}=CD\underline{L}=ABCD\underline{M},  \label{ej2}
\end{equation}
is obtained by eliminating from $H_{4}$ the columns $\bar{d}=\left\{
AB,AC,BC\right\} ,$ which form a set of rank 2. The WLP of (\ref{ej2}) is 
\begin{equation}
(0,0,16,39,48,48,48,39,16,0,0,1).  \label{wlpej2}
\end{equation}
If instead of eliminating the columns indicated above we choose $%
\{AB,AC,ABCD\},$ which has rank 3, the fraction obtained has as WLP 
\[
(0,0,17,38,44,52,54,33,12,4,1) 
\]
that has more words of length 3 and is therefore worse from an aberration
point of view.

When the number of columns eliminated is $h=2^{v}-1,$ where $v$ is any
integer, the choice of a set $\bar{d}$ with minimum rank ensures that the
resulting fraction has minimum aberration. This is due to the fact that all
sets of minimum rank and $2^{v}-1$ elements are isomorphic, Tang and Wu
(1996).

When $k>n/2$, the number of generators needed to define a design is larger
than the number of odd combinations of the basic factors. We will show that
the minimum aberration design is obtained by adding $r=k-n/2$ additional
columns to $O_{m}$. In this case, $q$ out of the $p$ generators needed are
the same as those used in the saturated design (\ref{saturado}). Hence, we
just have to find the additional $r=p-q$ generators among the $2^{m-1}-1$
free columns of $H_{m}$, which correspond to the even products of the basic
ones $B_{1},B_{2},...,B_{m}.$ We shall denote these columns as $\psi _{i}$,
to distinguish them from the odd ones $\varphi _{j}$; the set they form will
be denoted as 
\begin{equation}
E_{m}=\left\{ \psi _{1},\psi _{2},...,\psi _{2^{m-1}-1}\right\} .
\label{E_m}
\end{equation}
$E_{m}\subset H_{m}$ and $\left\{ I\right\} \cup E_{m}$ forms an abelian
subgroup of $\left\{ I\right\} \cup H_{m}.$ The sets $E_{m}$ and $H_{m-1}$
are isomorphic. The problem of choosing $r$ columns from within the $%
2^{m-1}-1$ even columns $(E_{m})$ is the same as the initial problem
(choosing\emph{\ }$k$ from $H_{m}$), though from a set of lower dimension $%
H_{m-1}$. Therefore, the problem can be solved through a iterative procedure.

\begin{theorem}
Let $d\subset H_{m}$ be a $2^{k-p}$ design with $n$ runs and $k>n/2.$ A
necessary condition for $d$ to be of MA is $O_{m}\subset d.$
\end{theorem}

\emph{Proof.} Let $H_{m}=d\cup \bar{d},$ where the number of columns of $%
\bar{d}$ is $h=(2^{m}-1)-k.$ Since $k>n/2,$ then $h<n/2=2^{m-1}$ and, by
theorem 1, the rank $v$ of $\bar{d}$ should verify $v\leq m-1$ for $d$ to be
of MA. The set $E_{m}$ is isomorphic to any set with all its $2^{m-1}-1$
elements generated by $m-1$ independent columns. Let $\bar{d}^{\prime }$ be
the set isomorphic to $\bar{d}$ such that $\bar{d}^{\prime }\subset E_{m}.$
The complementary of $\bar{d}^{\prime }$ contains all the odd combinations
of the basic columns and is isomorphic to $d\mathbf{.}$ $\square $

Another equivalent statement for the above theorem in terms of the defining
relation is the following: any minimum aberration $2^{k-p}$ design $d$ with $%
n=2^{k-p}$ runs and $k>n/2$ can be defined by 
\begin{eqnarray}
I &=&\varphi _{1}B_{m+1}=\varphi _{2}B_{m+2}=\cdots =\varphi _{q}B_{m+q}= 
\nonumber \\
&=&\psi _{1}^{\prime }B_{m+q+1}=\psi _{2}^{\prime }B_{m+q+2}=\cdots =\psi
_{r}^{\prime }B_{m+q+r}  \label{teorema_1}
\end{eqnarray}
where $\varphi _{1},\varphi _{2},...,\varphi _{q}$ are all $q=n/2-m$ columns
obtained as odd products of three or more of the basic columns $%
B_{1},B_{2},...,B_{m}$ and $\psi _{1}^{\prime },\psi _{2}^{\prime },...,\psi
_{r}^{\prime }$ are even products of the basic set.

\section{Saturated designs of resolution IV $(O_{m})$}

The strategy we presented in the next section to build resolution $III$
designs takes $O_{m}\subset H_{m},$ the saturated designs $2^{l-q}$ of
resolution IV with generators (\ref{saturado}), as the starting point. For
these designs $l=2^{m-1}=q+m.$ The smallest of these designs is $2^{4-1}$,
followed by $2^{8-4}$, $2^{16-11}$, $2^{32-26}$ and so forth.

The $2^{4-1}$ design defined by the columns $O_{3}=\{A,B,C,ABC\}$ has the
generating equation $I=ABC$\underline{$D$} and WLP $(0,0,0,1)$. The next
design is $O_{4},$ a $2^{8-4}$ fraction with generators 
\begin{equation}
I=ABC\underline{E}=ABD\underline{F}=ACD\underline{G}=BCD\underline{H}
\label{d8_4}
\end{equation}
and WLP: $(0,0,0,14,0,0,0,1).$

The defining relation is one of the $2^{m}$ \emph{alias chains} which appear
in a design $d\subset H_{m}$. The others correspond to the $2^{m}-1$ columns
or elements of $H_{m}$. Each chain includes $2^{q}$ confounded effects. The
set of these confusion chains is called the \emph{confusion structure}. The
confusion structure of $O_{m}$ designs with generators given in (\ref
{saturado}), presents certain special properties of interest. The $q$
generators of these fractions are even and therefore the $2^{q}$ words of
the defining relation are even. The alias chains associated to a main effect
are obtained by multiplying the generating equation by the single letter
associated to the factor; the result is a chain with $2^{q}$ odd effects.
Let $c_{2r+1,j}$ be the number of effects of size $2r+1$ in the alias chains
of $B_{j}$, it can be seen that for all $r$, 
\begin{equation}
c_{2r+1,1}=c_{2r+1,2}=\cdots =c_{2r+1,l}=\binom{l}{2r+1}/l.  \label{c_2r}
\end{equation}

On the other hand, the $(l-1)$ columns of $H_{m}$ not included in the set $%
O_{m}$ are the product of an even number of basic columns $\psi _{j}\in
E_{m} $. For the same reason as before, the alias chains associated to these
columns are formed by $2^{q}$ effects which are even. If $b_{2r,j}$ is the
number of effects of size $2r$ in the alias chains corresponding to the even
effect $\psi _{j}\in E_{m},$ then for any size $2r,$ 
\begin{equation}
b_{2r,1}=b_{2r,2}=\cdots =b_{2r,l-1}  \label{b_2r}
\end{equation}
These chains, which we shall call \emph{even chains}, play a fundamental
role in the computation of the word-length pattern for designs of resolution 
$III$ with generators given in (\ref{teorema_1}). Given $a_{2r}$, the value $%
b_{2r}\equiv b_{2r,j}$ is obtained through the following argument: the total
number of effects of size $2r$ which appear in the whole confusion structure
of the fraction is $\binom{l}{2r},$ $a_{2r}$ of which are words of the
defining relation and the remaining effects are distributed in the $(l-1)$
even alias chains, i.e. 
\begin{equation}
a_{2r}+(l-1)b_{2r}=\binom{l}{2r},\quad r=1,2,...,l/2  \label{c2i}
\end{equation}
Since $l=2^{m-1},$ the values $a_{2r}$ and $b_{2r}$ only depend on $m$ and
are easily obtained from the above expression and the following theorem.

\begin{theorem}
\label{t1} The \emph{word length pattern,} 
\[
W=(a_{1},a_{2},...,a_{l}), 
\]
of $O_{m},$ the saturated design $2^{l-q}$ of resolution IV, where $%
l=2^{m-1} $ and $q=2^{m-1}-m$, is obtained by means of the following
recurrence laws: 
\begin{eqnarray*}
a_{2r+1}=0,\quad r\geq 0; \\
a_{2}\quad =0, \\
a_{2r+2}=\dfrac{1}{2r+2}\left[ \dbinom{l}{2r+1}-(l-2r)a_{2r}\right] ,\qquad
(r\geq 1).
\end{eqnarray*}
\end{theorem}

\emph{Proof}. The design generators are even words and the product of any
subset of them is an even word, hence in the defining relation there exists
no even word and $a_{2r+1}=0$ for all $r.$

Since there are $a_{2r}$ words of size $2r$, the total number of letters
which appear in words of this size is $2ra_{2r}.$ If $\alpha _{2r,j}$ is the
number of times in which letter $B_{j}$ appears in words of size $2r$, we
have 
\begin{equation}
\sum_{j=1}^{l}\alpha _{2r,j}=2ra_{2r},  \label{adr}
\end{equation}
and according to (\ref{c_2r}) and (\ref{b_2r}), $\alpha _{2r,1}=\alpha
_{2r,2}=\cdots =\alpha _{2r,l},$ for all $r$, therefore 
\[
\alpha _{2r,j}=\frac{2r}{l}a_{2r}. 
\]
and the number of words of the same size $2r$ which do not contain the
letter $B_{j}$ is $a_{2r}-\alpha _{2r,j}=(1-2r/l)a_{2r}$.

To obtain the confusion chain corresponding to $B_{j}$, we multiply each
word of the generating relation by that factor. All the effects obtained
will be odd. The effects of size $2r+1$ in this chain correspond to either
words of size $2r$ in the defining relation which do not include factor $%
B_{j}$ or words of size $2r+2$ in the defining relation which do include
factor $B_{j}$. If $c_{2r+1,j}$ is the number of effects of size $2r+1$ in
the chain, we have 
\begin{eqnarray*}
c_{2r+1,j} &=&(a_{2r}-\alpha _{2r,j})+\alpha _{2r+2,j} \\
&=&\frac{l-2r}{l}a_{2r}+\frac{2r+2}{l}a_{2r+2}.
\end{eqnarray*}
All these odd effects can be found in the confusion chains of the main
effects. The total number of effects of size $2r+1$ is $%
\sum_{j=1}^{l}c_{2r+1,j}=\dbinom{l}{2r+1}$ and therefore 
\[
\left( l-2r\right) a_{2r}+(2r+2)a_{2r+2}=\dbinom{l}{2r+1}. 
\]
Solving for $a_{2r+2}$ we obtain the recurrence law. Since the design is of
resolution IV, the initial value of the law is $a_{2}=0.\square $

\section{Word-length pattern}

We will show a simple way of obtaining the WLP of a $2^{k-p}$ design with
generators (\ref{teorema_1}). For this purpose we shall now introduce the
concept of \emph{word-length pattern polynomial (WLPP).}

\begin{definition}
Let $d\subset H_{m}$ be a $2^{k-p}$ design with WLP, $%
W=(a_{1},a_{2},...,a_{k}).$ The polynomial $P_{d}$ associated to this design
is given by 
\[
P_{d}(u)=1+a_{1}u+a_{2}u^{2}+\cdots +a_{k}u^{k}. 
\]
We shall refer to it as the \emph{word length pattern polynomial} (WLPP).
\end{definition}

The above definition is valid for the case of full designs where $P_{d}(u)=1$%
.

For simplicity we shall denote by $P_{m}\equiv P_{O_{m}}$ the polynomial
associated to $O_{m}$, whose coefficients $a_{2r\text{ }}$ can be obtained
by applying Theorem \ref{t1}. Let us call 
\begin{equation}
Q_{m}(u)=\sum\limits_{r=1}^{2^{m-1}}b_{2r}u^{2r}  \label{Q_k}
\end{equation}
the polynomial with coefficients $b_{2r}\equiv b_{2r,j}$ defined in the
preceding section. The equation (\ref{c2i}) can be rewritten as 
\begin{equation}
P_{m}(u)+(2^{m-1}-1)Q_{m}(u)=\sum_{r=0}^{2^{m-2}}\binom{2^{m-1}}{2r}u^{2r}.
\label{S_k}
\end{equation}
We will refer to $Q_{m}$ as the \emph{effect-length pattern polynomial} for
the even alias chain.

Let $d\subset H_{m}$ be the $2^{k-p}$ design defined as $d=O_{m}\cup e,$
where $O_{m}$ is the design (\ref{saturado}) and $e\subset E_{m}$ with $%
r=k-n/2$ columns. Consider the \emph{design} $e$, a design formed with
columns from $E_{m}.$ The set $E_{m}$ is isomorphic to $H_{m-1},$ therefore
any design $e$ has a maximum of $m-1$ independent columns and is isomorphic
to a design of $H_{m-1}$. The word length pattern of $e$ can be obtained
taking into account the definition of word given in (\ref{word}) or,
alternatively, it could be find the isomorphic set of $H_{m-1}$ and write it
as in (\ref{d}). For example, the design with generators (\ref{ej2})
contains both $O_{4}$ and the subset 
\[
e=\left\{ AD,BD,CD,ABCD\right\} 
\]
of 
\[
E_{4}=\left\{ AB,AC,AD,BC,BD,CD,ABCD\right\} . 
\]
$E_{4}$ is formed by all the possible even combinations of the basic columns 
$A,B,C,D$. It can be seen that the rank of $E_{4}$ is 3. Using the
additional factors in (\ref{ej2}) to define the generators $J=AD,$ $K=BD,$ $%
L=CD$ and $M=ABCD,$ the set $e$ is isomorphic to $\{J,K,L,JKL\},$ which
defines a $2^{4-1}$ design of $H_{3}$. The defining relation of $e$ is $%
I=JKLM$ and hence its WLPP is $P_{e}(u)=1+u^{4}.$

\begin{theorem}
\label{th3} Let $d$ be a $2^{k-p}$ design defined by $d=O_{m}\cup e,$ where $%
O_{m}\subset H_{m}$ is the saturated design of resolution IV and $e\subset
E_{m}$. The (WLPP) of $d$ is given by 
\begin{equation}
P_{d}(u)=(1+u)^{r}Q_{m}(u)+P_{e}(u)[P_{m}(u)-Q_{m}(u)],  \label{teo_3}
\end{equation}
where $r=k-n/2,$ $P_{m}$ and $P_{e}$ are the WLPP of $O_{m}$ and $e,$
respectively, and 
\[
Q_{m}(u)=\frac{1}{2^{m-1}-1}\left[ \sum_{i=0}^{2^{m-2}}\binom{2^{m-1}}{2i}%
u^{2i}-P_{m}(u)\right] . 
\]
\end{theorem}

\emph{Proof}. Suppose that the design $d$ has only one more column than $%
O_{m}.$ In this case, the generators of $d$ are the $q=n/2-m$ generators of $%
O_{m}$ given in (\ref{saturado}), plus $B_{m+q+1}\psi _{1}^{\prime },$ where 
$\psi _{1}^{\prime }$ has an even number of basic factors. The defining
relation has $2^{q+1}$ words, the first $2^{q}$ correspond to the design $%
O_{m}$ which has WLP polynomial $P_{m}$. The remaining $2^{q}$ words are the
result of multiplying the $2^{q}$ words of the defining relation of $O_{m}$
by $B_{m+q+1}\psi _{1}^{\prime }.$ When multiplying the defining relation of 
$O_{m}$ by $B_{m+q+1}\psi _{1}^{\prime }$ we obtain the alias chain
associated to $\psi _{1}^{\prime }$ plus the new letter $B_{m+q+1}$ in each
effect. The effects of size $2i$ in the alias chain of $\psi _{1}^{\prime }$
become words of size $2i+1$ in the defining relation of $d.$ Since the
effect length pattern polynomial of the alias chain of $\psi _{1}^{\prime }$
is $Q_{m}(u),$ it will be $uQ_{m}(u)$ for the chain $B_{m+q+1}\psi
_{1}^{\prime }$ so that 
\[
P_{d}(u)=P_{m}(u)+uQ_{m}(u). 
\]

Assume now that $d$ is given by $O_{m}$ plus two additional columns, i.e.,
we add $B_{m+q+1}\psi _{1}^{\prime }$ and $B_{m+q+2}\psi _{2}^{\prime }$ to
the generators of $O_{m}.$ The defining relation of $d$ is formed by the
defining relation of $O_{m}$ plus the words resulting from multiplying the
defining relation of $O_{m}$ by $B_{m+q+1}\psi _{1}^{\prime },B_{m+q+2}\psi
_{2}^{\prime }$ and those corresponding to the product of both, $%
B_{m+q+1}B_{m+q+2}\psi _{1}^{\prime }\psi _{2}^{\prime }.$ The product $\psi
_{1}^{\prime }\psi _{2}^{\prime }$ is an element of $E_{m}$ different from $%
\psi _{1}^{\prime }$ and $\psi _{2}^{\prime }.$ The contributions of each
new generator to the generating equation of $d$ will be: $uQ_{m}(u)$,
corresponding to $B_{m+q+1}\psi _{1}^{\prime }$; the same value $uQ_{m}(u)$,
corresponding to $B_{m+q+2}\psi _{2}^{\prime }$ and $u^{2}Q_{m}(u)$ for $%
B_{m+q+1}B_{m+q+2}\psi _{1}^{\prime }\psi _{2}^{\prime }.$ The term $u^{2}$
is due to the fact that all the effects of the alias chain $\psi
_{1}^{\prime }\psi _{2}^{\prime }$ will now contain the two new letters $%
B_{m+q+1}B_{m+q+2}.$ The WLP polynomial of $d$ is 
\[
P_{d}(u)=P_{m}(u)+2uQ_{m}(u)+u^{2}Q_{m}(u) 
\]
adding and substracting $Q_{m}(u)$ to the right hand side, we obtain 
\[
P_{d}(u)=(1+u)^{2}Q_{m}(u)+[P_{m}(u)-Q_{m}(u)]. 
\]

The above result is generalized when $d=O_{m}\cup e,$ and the elements of $%
e=\left\{ \psi _{1}^{\prime },\psi _{2}^{\prime },...,\psi _{r}^{\prime
}\right\} $ are independent $(r\leq m-1).$ The independence of the elements
of $e$ ensures that any product of a subset of them provides a different
element of $E_{m}.$ With $r$ generators one can get $\binom{r}{i}$
combinations of $i$ of them that will incorporate $\binom{r}{i}u^{i}Q_{m}(u)$
to the polynomial $P_{d}(u),$%
\begin{eqnarray}
P_{d}(u) &=&P_{m}(u)+\sum_{i=1}^{r}\binom{r}{i}u^{i}Q_{m}(u)  \label{P_di} \\
&=&(1+u)^{r}Q_{m}(u)+[P_{m}(u)-Q_{m}(u)]  \nonumber
\end{eqnarray}

When the elements of $e=\left\{ \psi _{1}^{\prime },\psi _{2}^{\prime
},...,\psi _{r}^{\prime }\right\} $ are not independent (this is the case
where $r>m-1$), the product of some subset of them is equal to $I$ instead
of some element of $E_{m}.$ Thus, if 
\[
\psi _{j_{1}}^{\prime }\psi _{j_{2}}^{\prime }\cdots \psi _{j_{i}}^{\prime
}=I, 
\]
the product of the generators $B_{m+q+j_{1}}\psi _{j_{1}}^{\prime }$, $%
B_{m+q+j_{2}}\psi _{j_{2}}^{\prime },$ \ldots , $B_{m+q+j_{i}}\psi
_{j_{i}}^{\prime }$ is 
\[
B_{m+q+j_{1}}B_{m+q+j_{2}}\cdots B_{m+q+j_{i}} 
\]
and its contribution to $P_{d}(u)$ is $u^{i}P_{m}(u)$ instead of $%
u^{i}Q_{m}(u).$ Let $a_{i}^{\prime }$ be the number of subsets of $i$
elements of $e$ the product of which is equal to $I.$ The polynomial $%
P_{d}(u)$ is 
\begin{equation}
P_{d}(u)=P_{m}(u)+\sum_{i=1}^{r}\binom{r}{i}u^{i}Q_{m}(u)+%
\sum_{i=1}^{r}a_{i}^{\prime }u^{i}[P_{m}(u)-Q_{m}(u)],  \label{P_dd1}
\end{equation}
where the two first terms from (\ref{P_di}) correspond to the counting
process if the $r$ elements of $e$ are independent, and the term $%
\sum_{i=1}^{r}a_{i}^{\prime }u^{i}[P_{m}(u)-Q_{m}(u)]$ corrects this value
by taking into account the dependencies on $e.$ Once again, adding and
substracting $Q_{m}(u)$ to the right-hand side of (\ref{P_dd1}) and
regrouping terms one gets 
\[
P_{d}(u)=(1+u)^{r}Q_{m}(u)+(1+\sum_{i=1}^{r}a_{i}^{\prime
}u^{i})[P_{m}(u)-Q_{m}(u)]. 
\]
The polynomial $P_{e}(u)=(1+\sum_{i=1}^{r}a_{i}^{\prime }u^{i})$ is by
definition the WLPP of design $e.$ This last result is the expression (\ref
{teo_3}) which we wanted to prove.$\square $

\emph{Example 2}. The WLP for the $2^{12-8}$ minimum aberration fraction
with generators (\ref{ej2}) is obtained by applying (\ref{teo_3}), 
\[
P_{d}(u)=(1+u)^{4}Q_{4}(u)+P_{e}(u)\left[ P_{4}(u)-Q_{4}(u)\right] 
\]
where $P_{4}(u)=1+14u^{4}+u^{8}$ is the polynomial corresponding to the
fraction (\ref{d8_4}), and $Q_{4}(u)=4u^{2}+8u^{4}+4u^{6}$ is obtained from (%
\ref{S_k}) and (\ref{Q_k}). Also, $P_{e}(u)=1+u^{4}$ is the polynomial of
the $2^{4-1}$ minimum aberration design ($I=JKLM)$. By substitution in $%
P_{d} $ we obtain: 
\[
P_{d}(u)=1+16u^{3}+39u^{4}+48u^{5}+48u^{6}+48u^{7}+39u^{8}+16u^{9}+u^{12}. 
\]
The coefficients of $P_{d}$ are the elements of the word length pattern (\ref
{wlpej2}).

Using theorem \ref{th3}, it is easy to prove that the goodness of the design 
$d$ in terms of its aberration depends entirely on $e.$

\begin{theorem}
Let $d$ be a design defined by $d=O_{m}\cup e,$ where $O_{m}\subset H_{m}$
is the saturated design of resolution IV and $e\subset E_{m}$. The design $d$
is MA if and only if $e$ is isomorphic to an MA design $e^{\prime }\subset
H_{m-1}$.
\end{theorem}

\emph{Proof. }Let $d_{1}\subset H_{m}$ and $d_{2}\subset H_{m}$ be two
different $2^{k-p}$ fractions. The design $d_{1}$ is of smaller aberration
than $d_{2}$ if and only if (by definition of WLPP) $%
P_{d_{1}}(u)-P_{d_{2}}(u)$ has a negative coefficient for the term of lowest
order. If $d_{1}=O_{m}\cup e_{1}$ and $d_{2}=O_{m}\cup e_{2}$, then by (\ref
{teo_3}) we have 
\[
P_{d_{1}}(u)-P_{d_{2}}(u)=[P_{e_{1}}(u)-P_{e_{2}}(u)][P_{m}(u)-Q_{m}(u)] 
\]
and since $P_{m}(u)-Q_{m}(u)=1-b_{2}u^{2}+\cdots ,$ if the coefficient of
the lowest order term of $P_{e_{1}}(u)-P_{e_{2}}(u)$ is negative, the
corresponding one for $P_{d_{1}}(u)-P_{d_{2}}(u)$ will also be negative.
Therefore, if $e_{1}$ is isomorphic to an MA design of $H_{m-1}$, $d_{1}$
will also be MA, and vice versa. $\square $

This theorem is useful to justify that the design $d$ with generators (\ref
{ej2}) is MA. This follows from the fact that this design is the union of $%
O_{4}$ and $e=\left\{ AD,BD,CD,ABCD\right\} ,$ and the subset $e$ is
isomorphic to $\{J,K,L,JKL\},$ which represents the $2^{4-1}$ MA design.

Knowing the $2^{4-1}$ MA $(m=3)$ design it is possible to obtain the $%
2^{12-8}$ MA $(m=4)$ design. This design in turn allows the determination of
the $2^{28-23}$ MA $(m=5)$ design, and so forth. Given the MA design $%
d\subset H_{m}$ with $k$ factors, one can immediately obtain the MA design $%
d^{\prime }\subset H_{m+1}$ with $2^{m}+k$ factors.

The $d\subset H_{m}$ MA design with a number of factors $k>n/2,$ is obtained
as $d=O_{m}\cup e,$ where $O_{m}$ is the (\ref{saturado}) design and $%
e\subset E_{m}$ with $r=k-n/2$ columns. To determine $e$ one should obtain
the MA with $r$ factors in $H_{m-1}.$ Let $d^{\prime }\subset H_{m-1}$ be
the MA design with $r$ factors defined from the basic factors $B_{1}^{\prime
},B_{2}^{\prime },...,B_{m-1}^{\prime }.$ The set $e$ is defined as the
transform of $d^{\prime }$ by the isomorphism $\mathcal{T}$ from $H_{m-1}$
to $E_{m}$ defined by $\mathcal{T}(B_{j}^{\prime })=B_{j}B_{m}$ for $%
j=1,2,...,m-1.$ When $r>n/4$ the method can be applied iteratively.

\emph{Example 3}. Let us illustrate, for instance, how to obtain the $%
2^{28-23}$ MA design. This design is included in $H_{5},$ and, following
theorem 2, its first 16 columns correspond to all the odd combinations $%
(O_{5})$ of the basic columns $A,B,C,D,E$. The remaining $12$ columns are a
subset of the even combinations of the basic ones $(E_{5})$ isomorphic to
the $2^{12-8}$ MA design included in $H_{4}$. The $2^{12-8}$ MA design
(verifying $12>n/2=8$) should contain the 8 odd combinations of the 4
independent columns plus 4 other combinations of them. We take $%
e_{0}=\left\{ AE,BE,CE,DE\right\} $ as the independent columns of $E_{5}$
and build with them the 8 odd combinations, 
\[
e_{1}=\{AE,BE,CE,DE,ABCE,ABDE,ACDE,BCDE\}. 
\]
The 7 even combinations of $e_{0}$ form the set $\{AB,AC,AD,BC,BD,CD,ABCD\}$
isomorphic to $H_{3}.$ From it we take 
\[
e_{2}=\{AD,BD,CD,ABCD\} 
\]
which is isomorphic to the $2^{4-1}$ MA design of $H_{3}$. The $2^{28-23}$
MA design is $d=O_{5}\cup e_{1}\cup e_{2}.$

The case $2^{28-23}$ has been chosen to illustrate the proposed method. A
simpler way to build the design in this case is to eliminate $3$ columns out
of the 31 which form $H_{5}.$ By results of section 2, it suffices to choose
3 minimum rank columns of $H_{5}$. In this case, $d$ is missing just the
columns $\bar{d}=\left\{ AB,AC,BC\right\} $ to complete $H_{5}$, and the
rank of $\bar{d}$ is the minimum possible.

\section{Conclusions}

The selection of a $2^{k-p}$ regular fraction implies choosing $k$ columns
among the $n=2^{m}$ available in a full factorial design with $m=k-p$
factors. When $k\geq n/2,$ it can be shown that the minimum aberration
design contains $n/2$ columns obtained as a product of an even subset of $m$
independent columns. The choice of the other $k-n/2$ columns requires
solving the same problem but now in a smaller set, the full factorial design
with $m-1$ factors. A iterative procedure has been derived from this
property that allows building screening designs ($k>n/2)$ of minimum
aberration for large $k.$

The comparison of designs in terms of aberration requires the computation of
word-length patterns. When $k$ and $p$ are large the direct counting process
is very cumbersome. In this article we introduce a simple relationship,
based on products of polynomials, for computing the word-length pattern for
these designs. This relation provides an interesting correspondence between
designs of different number of runs.

\section{Appendix: proof of theorem 1.}

We choose $\bar{d}=\left\{ \bar{\xi}_{1}^{\prime },\bar{\xi}_{2}^{\prime
},...,\bar{\xi}_{h}^{\prime }\right\} $ so that the complementary design $d$
has the minimum number of words of size 3 in its defining relation. In the
complete defining relation of $H_{m},$ any factor (letter) appears in $%
(n-2)/2$ words of size 3, therefore, if we remove $h$ columns, the number of
words of size 3 that we eliminate is 
\begin{equation}
a(\overline{d})=(\frac{n-2}{2})h-a_{3}^{\prime }(\overline{d})-2a_{3}(%
\overline{d}),  \label{coincidencias}
\end{equation}
where $a_{3}^{\prime }(\overline{d})$ and $a_{3}(\overline{d})$ are,
respectively, the number of these words which contain two and three
eliminated factors. Moreover, in the complete generating equation of $H_{m}$
any combinations of the two factors appear once in words of size 3, thus 
\[
\binom{h}{2}=a_{3}^{\prime }(\overline{d})+3a_{3}(\overline{d}). 
\]
Substituting in (\ref{coincidencias}) we find that the number of words of
size 3 that have been eliminated is 
\[
a(\overline{d})=\frac{h(n-h-1)}{2}+a_{3}(\overline{d}). 
\]
and this number will be maximum whenever $a_{3}(\overline{d})$ is also as
large as possible.

The term $a_{3}(\overline{d})$ denotes the number of relations or words 
\begin{equation}
\bar{\xi}_{t}^{\prime }\bar{\xi}_{u}^{\prime }\bar{\xi}_{v}^{\prime }=I,%
\text{ with }\bar{\xi}_{t}^{\prime },\bar{\xi}_{u}^{\prime },\bar{\xi}%
_{v}^{\prime }\in \overline{d}.  \label{trio}
\end{equation}
We shall now see that if $a_{3}(\overline{d})$ is maximum, then $\overline{d}
$ is generated by the smallest number of independent columns $v_{h}\in N,$
with $h<2^{v_{h}}$ (note that it is not possible to generate $h$ different
columns using less than $v_{h}$ independent columns).

Suppose that the number of independent elements of $\bar{d}$ is $v>v_{h},$
we shall now see that it is then possible to define a new set with $h$
columns which provides a larger number of words of size 3, which contradicts
the statement that $a_{3}(\overline{d})$ is maximum.

Let $J_{v}$ be the set formed by $2^{v}-1$ different columns generated by
the $v$ independent columns of $\bar{d}$. Out of the $v$ independent columns
we choose $v-1$, which we will call $\psi _{1},\psi _{2},...,\psi _{v-1}$.
Let $J_{v-1}\subset J_{v}$ be the subset generated by these columns and $%
\bar{J}_{v-1}$ the complementary set, such that $\bar{J}_{v-1}=J_{v}%
\backslash J_{v-1}.$ We divide the set $\overline{d}_{0}\equiv \overline{d}$
into two parts, 
\[
e_{0}=\overline{d}_{0}\cap J_{v-1}=\{\psi _{1},\psi _{2},...,\psi _{s}\} 
\]
and 
\[
f_{0}=\overline{d}_{0}\cap \bar{J}_{v-1}=\{\varphi _{1},\varphi
_{2},...,\varphi _{h-s}\}. 
\]
For this proof we have used on purpose the notation $\psi _{i},$ $\varphi
_{j}$ for the columns of $J_{v},J_{v-1},$ respectively, which, without loss
of generality, can be considered as even and odd combinations of the basic
factors. As will be seen below, this notation helps to understand the
following development, if we bear in mind properties such as that the
product of two even combinations is also even and so forth.

The elements of $f_{0}$ $\subset J_{v}$ are obtained as a product of the
generators of $J_{v-1}$ $\psi _{1},\psi _{2},...,\psi _{v-1}$ and any other
term of $f_{0},$ for example $\varphi _{1},$ thus 
\begin{equation}
\varphi _{t}=\varphi _{1}\varpi _{t},\quad \mathrm{con\quad }\varpi _{t}\in
J_{v-1},\text{\textrm{para }}t=2,...,h-s.  \label{phis}
\end{equation}
With the new notation, the relations (\ref{trio}) can be divided into two
classes: 
\begin{equation}
(\mathrm{A):\ }\psi _{i}\psi _{j}\psi _{k}=I,\quad \psi _{i},\psi _{j},\psi
_{k}\in e_{0}  \label{tipo A}
\end{equation}
and 
\begin{equation}
(\mathrm{B):\ }\varphi _{i}\varphi _{j}\psi _{k}=I\quad \varphi _{i},\varphi
_{j}\in f_{0},\quad \psi _{k}\in e_{0}.  \label{tipo B}
\end{equation}
(Note that taking into account the even/odd features no other combinations
are possible).

We now replace $\varphi _{t}\in f_{0}$ by $\psi _{t}^{\prime }\in J_{v-1}$ ($%
\psi _{t}^{\prime }\notin e_{0})$ as follows: we choose two different $\psi
_{a},\psi _{b}\in e_{0}$ such that $\psi _{1}^{\prime }=\psi _{a}\psi
_{b}\notin e_{0},$ (since $s<2^{v-1}$ this choice is always possible). In (%
\ref{phis}), $\varphi _{1}$ is replaced by $\psi _{1}^{\prime }$ and we
obtain the columns $\psi _{t}^{\prime }\in J_{v-1},$ such that 
\[
\psi _{t}^{\prime }=\psi _{1}^{\prime }\varpi _{t}\quad ,t=2,3,...,h-s. 
\]
We define $e_{0}^{^{\prime }}=\{\psi _{t}^{\prime }:\psi _{t}^{\prime }=\psi
_{1}^{\prime }\varpi _{t}\notin e_{0}\},$ and the new sets $e_{1}=e_{0}\cup
e_{0}^{\prime }$ $\cup \left\{ \psi _{1}^{\prime }\right\} $ and 
\[
f_{1}=\{\varphi _{t}=\varphi _{1}\varpi _{t}\in f_{0}:\psi _{t}^{\prime
}=\psi _{1}^{\prime }\varpi _{t}\in e_{0},\quad t=2,3,...,h-s\}. 
\]
If $\overline{d}_{1}=e_{1}\cup f_{1},$ we now show that $a_{3}(\overline{d}%
_{1})\geq a_{3}(\overline{d}_{0})$. Since $e_{0}\subset e_{1},$ the other
words of type (\ref{tipo A}) initially formed can still be formed after the
substitution process. For the words of type (\ref{tipo B}), there are three
possible situations:

\begin{enumerate}
\item  $\varphi _{i},\varphi _{j}\in \overline{d}_{1}$, that is $\varphi
_{i},\varphi _{j}$ have not no been replaced in which case the word $\varphi
_{i}\varphi _{j}\psi _{k}=I$ is maintained.

\item  Both $\varphi _{i},\varphi _{j}\notin \overline{d}_{1},$ have been
replaced by $\psi _{i}^{\prime }=\psi _{1}\varpi _{i}$ and $\psi
_{j}^{\prime }=\psi _{1}\varpi _{j},$ respectively. In this case $\psi
_{i}^{\prime }\psi _{j}^{\prime }\psi _{k}=I$ will be a new word . Given
that $\varphi _{i}\varphi _{j}\psi _{k}=I,$ according to (\ref{phis}) $%
\varpi _{i}\varpi _{j}\psi _{k}=I,$ therefore $\left( \psi _{1}\varpi
_{i}\right) \left( \psi _{1}\varpi _{j}\right) \psi _{k}=\psi _{i}^{\prime
}\psi _{j}^{\prime }\psi _{k}=I.$

\item  $\varphi _{i}\notin \overline{d}_{1},\varphi _{j}\in \overline{d}%
_{1},\varphi _{i}$ have been replaced by $\psi _{i}^{\prime }=\psi
_{1}\varpi _{i},$ but $\varphi _{j}$ has not, because the substitute $\psi
_{j}^{\prime }$ already existed, i.e. $\psi _{j}^{\prime }=\psi _{1}\varpi
_{j}\in \overline{d}_{0}.$ Using the same reasoning as above, since $\varphi
_{i}\varphi _{j}\psi _{k}=I,$ then $\psi _{i}^{\prime }\psi _{j}^{\prime
}\psi _{k}=I$ will be a word that did not exist before (given that $\psi
_{i}^{\prime }$ is a new element).

In this case it is possible that two different words give equal transforms.
Suppose that $\varphi _{i},\varphi _{r},\varphi _{s}\in f_{0},\psi _{u},\psi
_{v}\in e_{0}$, that 
\begin{equation}
\varphi _{i}\varphi _{r}\psi _{u}=I\text{ and }\varphi _{i}\varphi _{s}\psi
_{v}=I,  \label{cruce}
\end{equation}
and $\varphi _{i}\notin \overline{d}_{1}$ has been replaced by $\psi
_{i}^{\prime }\in e_{1},$ while $\varphi _{r},\varphi _{s}\in f_{1}$ could
not be replaced. Moreover, it holds that $\psi _{u}^{\prime }=\psi
_{1}\varpi _{u}=\psi _{v}$ and $\psi _{v}^{\prime }=\psi _{1}\varpi
_{v}=\psi _{u}.$ In this situation, the two last words (\ref{cruce}) have
the same transforms $\psi _{i}^{\prime }\psi _{u}\psi _{v}=I.$ But in this
case the new word $\varphi _{r}\varphi _{s}\psi _{i}^{\prime }=I$ is
created, which did not exist since $\psi _{i}^{\prime }\notin \overline{d}%
_{0}$.

A third coincidence could appear, if $i=1$ and $\psi _{a}=\psi _{u}$ and $%
\psi _{b}=\psi _{v}.$ In this case the transform of (\ref{cruce}) coincides
with $\psi _{1}^{\prime }\psi _{a}\psi _{b}=I,$ the relation used to define $%
\psi _{1}^{\prime }.$
\end{enumerate}

The substitution process verifies that $a_{3}(\overline{d}_{0})\leq a_{3}(%
\overline{d}_{1}).$ If there is a triple coincidence, as described in point
3 of the analysis above, then $a_{3}(\overline{d}_{0})\leq a_{3}(\overline{d}%
_{1}),$ if not, $a_{3}(\overline{d}_{0})<a_{3}(\overline{d}_{1})$ (all words
have been transformed and additionally $\psi _{1}^{\prime }\psi _{a}\psi
_{b}=I).$

Repeating the process, a finite sequence of sets $\overline{d}_{0},\overline{%
d}_{1},\overline{d}_{2},...,\overline{d}_{M}$ ($M\leq h-q),$ is created,
which verify 
\begin{equation}
a_{3}(\overline{d}_{0})\leq a_{3}(\overline{d}_{1})\leq \cdots \leq a_{3}(%
\overline{d}_{M-1})<a_{3}(\overline{d}_{M})  \label{trios}
\end{equation}
where $\overline{d}_{M}\subset J_{v-1}$ ($f_{M}=\varnothing ).$ Moreover, in
the last step, $a_{3}(\overline{d}_{M-\mathbf{1}})<a_{3}(\overline{d}_{M})$
since the situation described in point 3 can not arise, since $%
f_{M}=\varnothing $ (all terms of type $\varphi _{t}$ have been replaced $).$

The final set $\overline{d}_{M}$ obtained has more words than $\overline{d}$%
, which contradicts the statement that $a_{3}(\overline{d})$ is maximum.
Therefore, the theorem is proved.$\square $

\section{\textbf{References}}

\parindent 0cm \parskip 5pt

Ankenman, B. E. (1999), ``Design of Experiments with Two- and Four-Level
Factors.'' \emph{Journal of Quality Technology, }\textbf{31, }363-375.

Bingham, D. and R. R. Sitter (1999), ``Minimum-Aberration Two-Level
Fractional Factorial Split-Plot Designs.'' \emph{Technometrics, }\textbf{41,}
62-70.

Box, G. E. P. and J. S. Hunter, (1961). ``The $2^{k-p}$ fractional factorial
designs.'' \emph{Technometrics.} \textbf{3,} 311-351

Box, G. E. P., W. G. Hunter, and J. S. Hunter, (1978). \emph{Statistics for
Experimenters}. Wiley, New York.

Chen, H. and A.S. Hedayat, (1996). ``$2^{n-l}$ designs with weak minimum
aberration.'' \emph{Annals of Statistics}. \textbf{24} 2536-2548.

Chen, J. (1992). ``Some results on $2^{n-k}$ fractional factorial designs
and search for minimum aberration designs.'' \emph{Annals of Statistics}. 
\textbf{20,} 2124-2141.

Chen, J., Sun, D.X., and Wu, C.F.J. (1993). ``A catalogue of two-level and
three-level fractional factorial designs with small runs.'' \emph{Internat.
Statist. Rev.} \textbf{61} 131-145.

Chen, J. and C. F. J. Wu, (1991). ``Some results on $s^{\prime n-k}$
fractional factorial designs with minimum aberration or optimal moments.'' 
\emph{Annals of Statistics}. \textbf{19,} 1028-1041.

Franklin, M. F. (1984). ``Constructing tables of minimum aberration $p^{n-m}$
designs.'' \emph{Technometrics}. \textbf{26 }225-232.

Fries, A. and W.G. Hunter (1980). ``Minimum aberration $2^{k-p}$ designs.'' 
\emph{Technometrics} \textbf{22} 601-608.

Huang, P., D. Chen, and J. O. Voelkel (1998), ``Minimum-Aberration Two-Level
Split-Plot Designs.'' \emph{Technometrics}, 40, 314-326.

Pu, K. (1989), ``Contributions to fractional factorial designs.'' Ph. D.
dissertation, Univ. Illinois, Chicago.

Tang, B. and Wu, C. F. J. (1996) ``Characterization of minimum aberration $%
2^{n-k}$ designs in terms of their complementary designs.'' \emph{Annals of
Statistics}. \textbf{24} 2549-2559.

Tiao, G.C., Bisgaard S., Hill, W.J., Pe\~{n}a, D. and Stigler, S.M. (eds.)
(2000), ``\textit{Box on Quality and Discovery}.'' John Wiley \& Sons, New
York.

\end{document}